\documentclass{amsart}
\usepackage{amsmath}
\usepackage{amssymb} 
\usepackage{enumerate}

\newcommand{\dis}{\displaystyle}

\newtheorem{theorem}{Theorem}
\newtheorem{proposition}[theorem]{Proposition}
\newtheorem{question}[theorem]{Question}

\begin{document}

\title{A note on genera of band--connected sums that are fibered}
\author[K. Miyazaki]{Katura Miyazaki}
\address{Faculty of Engineering, Tokyo Denki University, 5 Senju Asahi-cho, Adachi-ku, Tokyo 120--8551, 
Japan}
\email{miyazaki@cck.dendai.ac.jp}

\date{}


\begin{abstract}
We show that
if a fibered knot $K$ is expressed as a band--connected sum of
$K_1, \ldots, K_n$,
then each $K_i$ is fibered, and the genus of $K$ is greater
than or equal to that of the connected sum of $K_1,\ldots,K_n$.
\end{abstract}

\maketitle


A \textit{band} $b$ attached to a link is an embedded disk in $S^3$ intersecting
the link in two subarcs of $\partial b$.
Let $L =K_1\cup \cdots \cup K_n$ be an $n$-component split link,
and take mutually disjoint bands $b_1,\ldots, b_{n-1}$ attached to $L$ so that
$\dis L \cup \bigcup_{i=1}^{n-1} b_i$ is connected, and
the orientations of $L$ and $b_i$ are consistent.
The knot $\dis K = L\cup \bigcup_{i=1}^{n-1} \partial b_i
- \bigcup_{i=1}^{n-1} \mathrm{int}( L\cap b_i )$ is
a \textit{band--connected sum} of $K_1,\ldots, K_n$
along $b_1, \ldots, b_{n-1}$.
A natural question to ask is whether the band--connected sum $K$
is more ``complicated" than the connected sum $K_1\#\cdots \# K_n$.
In particular, is the genus of a band--connected sum greater than
that of the connected sum?
For band--connected sums of two knots,
by using the theory of sutured manifolds
Gabai \cite{Gabaiband}, Scharlemann \cite{Sch1}, 
and Kobayashi \cite{Koba92} obtain the following results.

\begin{theorem}
\label{band sum1}
Let $K$ be a band--connected sum of $K_1, K_2$ along a band $b$.
Then the following hold.
\begin{enumerate}[$(1)$]
\item \cite{Sch1, Gabaiband}
$g(K) \ge g(K_1) + g(K_2)$ holds,
where $g( \cdot)$ is the genus of a knot.
The equality holds if and only if $K$ bounds a Seifert surface that
is the union of minimal genus Seifert surfaces for $K_i$ $(i=1,2)$
 and $b$.

\item \cite{Koba92} If $K$ is a fibered knot, then $K_1$ and $K_2$ are also fibered.
The equality $g(K) =g(K_1) + g(K_2)$ holds if and only if
there is a sphere in $S^3 - K_1\cup K_2$ 
intersecting a core of $b$ in one point
and consequently $K \cong K_1 \# K_2$.
\end{enumerate}
\end{theorem}

For a fibered knot that is a band--connected sum of two or more knots,
we prove Theorem~2 below.

\begin{theorem}
\label{band sum2}
Let $K$ be a band--connected sum of $K_1,\cdots, K_n$.
If $K$ is a fibered knot, then the following hold.
\begin{enumerate}[$(1)$]
\item The knots $K_1,\ldots, K_n$ are all fibered.
\item $g(K) \ge g(K_1) + \cdots +g(K_n)$ holds.
The equality holds if and only if $K\cong K_1 \# \cdots \# K_n$.
\end{enumerate}
\end{theorem}

\begin{question}
If $K$ is not fibered, does the inequality in Theorem~$\ref{band sum2}(2)$ hold?
\end{question}

Let us say a fibered knot is \textit{tight} if its open book decomposition
of $S^3$ induces the tight contact structure on $S^3$.
By using Theorem~\ref{band sum1} and the Ozsv\'ath---Szab\'o concordance invariant $\tau$ \cite{OS03},
Baker and Motegi \cite{BaMo} show that
if a non--trivial band--connected sum of 
two knots is a tight fibered knot, 
then the band sum is the connected sum.
Applying Theorem~\ref{band sum2}
instead of Theorem~\ref{band sum1} in their proofs,
we can extend their results to band--connected sums of more than
two knots as follows.
(See Subsection~5.2 in \cite{BaMo}.)



\begin{proposition}
\label{tight band sum}
Let $K$ be a fibered knot that is a band--connected sum of
$K_1,\ldots,K_n$.
If $K\not\cong K_i$ for all $i$, then the following hold.
\begin{enumerate}[$(1)$]
\item If $K$ is a prime knot, then $K$ is not tight.
\item $K$ does not admit a lens space surgery or more generally 
an $L$-space surgery.
\end{enumerate}
\end{proposition}

Theorem~2 is proved by combining some results on ribbon concordance
and
the positive solution to Rapaport's conjecture on knot--like groups.

For knots $K_0$ and $K_1$, we say that
$K_1$ is \textit{ribbon concordant} to $K_0$
and write $K_1 \ge K_0$
if there is a concordance $C \subset S^3\times I$
between $K_i \subset S^3\times \{i\}$ $(i=0,1)$ on which
the projection $S^3\times I \to I$ is a Morse function
with no local maxima (Gordan \cite{Go81}).
Such $C$ is called a \textit{ribbon concordance} from $K_1$ to $K_0$.
\cite[Lemma~3.4]{Go81} implies that
if $K_1\ge K_0$ and both $K_i$ are fibered, then
$g(K_1) > g(K_0)$ or $K_1 \cong K_0$.
By \cite{Miya98} any band-connected sum of $K_1, \ldots, K_n$ is
ribbon concordant to the connected sum $K_1 \# \cdots \# K_n$.
Theorem~\ref{band sum2} thus follows from Proposition~\ref{Silver}.

\begin{proposition}
\label{Silver}
If $K_1 \ge K_0$ and $K_1$ is fibered, then $K_0$ is also fibered.
\end{proposition}

Silver \cite{Silver92} observed that
E. S. Rapaport's conjecture on knot--like groups \cite{Rapaport}
implies Proposition~\ref{Silver}.
After that the conjecture was proved by
Kochloukova \cite[Corollary~2]{Koch06} in 2006.
Proposition~\ref{Silver} and Theorem~\ref{band sum2}
is thus obtained.
For the sake of the readers,
here we present Rapaport's conjecture and Silver's argument.

Neuwirth \cite{Neuwirth} shows that for any knot group $G$
if its commutator subgroup $G'$ is finitely generated, then $G'$ is free.
Rapaport's conjecture \cite{Rapaport} states Neuwirth's theorem
holds more generally for any knot--like group;
a knot--like group is a finitely presented group that has
the abelianization $\mathbb{Z}$ and deficiency 1.
See \cite{Koch06} (and also \cite[Theorem~2.1]{Hill})
for the proof of the conjecture.

For a ribbon concordance $C \subset S^3\times I$ from $K_1$ to $K_0$,
set $G = \pi_1( S^3\times I -\mathrm{int} N(C) )$ and
$G_i = \pi_1( S^3\times \{i\} -\mathrm{int} N(K_i) )$, $i=0,1$.
The ribbon concordance group $G$ is a knot--like group because
it has the presentation
$\langle G_0, x_1, \ldots, x_n |\, r_1, \ldots, r_n \rangle$ \cite{Go81},
where $n$ is the number of local minima on $C$.
Also by \cite{Go81},
the inclusion maps $S^3\times \{ i \} -\mathrm{int}N(K_i)
\to S^3\times I -\mathrm{int}N(C)$ $(i=0,1)$ induce
an epimorphism $G_1 \to G$ and a monomorphism $G_0 \to G$.
We then have an epimorphism $G'_1 \to G'$ and
a monomorphism $G'_0 \to G'$.
Suppose $K_1$ is fibered.
Since $G_1'$ is finitely generated, $G'$ is also finitely generated.
By Rapaport's conjecture, $G'$ is a free group, so that
$G'_0$ is free.
It follows that $K_0$ is fibered.

\bigskip


\end{document}